\documentclass[a4paper,reqno]{amsart}
\usepackage{amssymb, amsmath, amscd}

\newcommand{\R}{\mathbb{R}}
\newcommand{\N}{\mathbb{N}}

\newcommand{\Z}{\mathbb{Z}}

\newtheorem{theorem}{Theorem}
\newtheorem{lemma}[theorem]{Lemma}

\newtheorem{remark} {Remark}

\newtheorem{conjecture}{Conjecture}

\def\text#1{\hbox{#1}}

\def\gronk{\vphantom{\vrule height 12pt}}
\def\operatorname#1{\hbox{#1}}
\def\proofx{\medbreak\noindent{\it Proof.}}
\def\qedbox{\hbox{$\rlap{$\sqcap$}\sqcup$}}
\def\eqref#1{Equation~(\ref{#1})}
\def\iint{\displaystyle\int\int}

\begin{document}
\title{Heat content asymptotics with singular data}
\author{M.van den Berg}
\address
{School of Mathematics, University of Bristol\\
University Walk, Bristol BS8 1TW\\United Kingdom}
\begin{email}{M.vandenBerg@bris.ac.uk}\end{email}
\author{P.Gilkey}
\address{Mathematics Department, University of Oregon\\
Eugene, OR 97403\\USA}
\begin{email}{gilkey@uoregon.edu}\end{email}
\begin{abstract}\noindent We study the asymptotic behaviour of the
heat content on a compact Riemannian manifold with boundary and with singular
specific heat and singular initial temperature distributions.
Assuming the existence of a complete asymptotic series we
determine the first three terms in that series. In addition to the general setting, the interval is studied in detail.
\par\noindent{\it Classification\/}: 58J32; 58J35; 35K20
\end{abstract}
\maketitle
\section*{Dedication} This paper is dedicated to Stewart Dowker. Stuart has made many contributions to the field of
spectral geometry and remains active in this area
\cite{D10,D11}. The second author has been honored to have been a collaborator with Stewart
\cite{DGK99,DKG01}. We hope this paper serves as a fitting tribute to our
colleague and friend.

\section{Introduction}
Let $M$ be a compact Riemannian manifold with smooth boundary $\partial M$, and let $\delta$ denote the geodesic
distance to the boundary. Let $\psi_1$ and $\psi_2$ be smooth functions on the interior of
$M$. Then $\psi_1$ will represent the initial temperature of $M$ and $\psi_2$ will represent the specific heat of $M$. Since $M$ is compact and $\partial M$ is smooth
the distance function is smooth near $\partial M$. We have to assume
that $\delta^{\alpha_1}\psi_1$ and $\delta^{\alpha_2}\psi_2$ are smooth on a closed collared neighbourhood of $\partial M$. The parameters $\alpha_1$ and $\alpha_2$
control the growth or decay of $\psi_1$ and $\psi_2$ near $\partial M$. Let
$D$ be an operator of Laplace type on
$M$. Impose Dirichlet boundary conditions to define the realization of
$D$. Let $dx$ be the Riemannian measure on $M$. Since $\partial M$ is smooth, the corresponding Dirichlet heat kernel
$p_M(x_1,x_2;t), x_1 \in M, x_2 \in M, t>0$ vanishes
linearly in $\delta(x_1)$ and $\delta(x_2)$ near $\partial M$. We suppose $\alpha_1<2$ and $\alpha_2<2$
to ensure convergence subsequently. Let $e^{-tD}$ be the fundamental solution of the heat equation for the
Dirichlet Laplacian. Then
$$u_1(\cdot,t):=e^{-tD}\psi_1=\int_Mp_M(\cdot,x;t)\psi_1(x)dx$$
 represents the temperature of the manifold for
$t>0$.  The {\it heat content} $Q$ is defined by
\begin{eqnarray*}
Q(\psi_1,\psi_2,D)(t):&=&\int_Mu_1(x;t)\cdot\psi_2(x)dx\\
&=&\int_M\int_Mp_M(x_1,x_2;t)\psi_1(x_1)\psi_2(x_2)dx_1dx_2\,.
\end{eqnarray*}

\begin{conjecture}\label{conj-1.1} Let $\alpha_1+\alpha_2\notin \Z$, $\alpha_1<2$, $\alpha_2<2$.
There is a complete asymptotic series as $t\downarrow0$
\begin{eqnarray*}
&&Q(\psi_1,\psi_2,D)(t)\sim\sum_{n=0}^\infty t^{n}\beta_n^M
+\sum_{j=0}^\infty t^{(1+j-\alpha_1-\alpha_2)/2}\beta_j^{\partial
M},
\end{eqnarray*}
where the $\beta_n^M, n=0,1,\cdots$ are regularized integrals of local invariants over
$M$  and where the $\beta_j^{\partial M}, j=0,1,\cdots$ are integrals of
local invariants over the boundary.
\end{conjecture}
\begin{remark}\rm It is convenient to let $\alpha_1$ and $\alpha_2$ be complex as we may then use analytic
continuation. For
$\Re(\alpha_1)<<0$ and
$\Re(\alpha_2)<<0$,
$$\beta_n^M=(-1)^n\frac1{n!}\int_MD^n\psi_1(x)\cdot\psi_2(x) dx\,.$$
The values of $\beta_n^M$ for more general values of $\alpha_1$ and $\alpha_2$ may then be obtained as
regularized integrals as discussed in \cite{vdB5}. We omit the technical details concerning the requisite regularizations in the
interests of brevity as they will play no role in our analysis.
\end{remark}

The heat content has obvious physical relevance and the invariants $\beta_j^{\partial
M}$, which reflect the asymptotic behaviour as $t\downarrow 0$, relate the geometry of $M$ to the underlying physical
properties of $M$. Much of the
previous work in the field has been devoted to the computation of the invariants
$\beta_j^{\partial M}$ in the smooth setting
($\alpha_1=0$, $\alpha_2=0$). They were originally studied for the scalar Laplacian with $\psi_1=\psi_2=1$ \cite{BD89,BlG94,BS90}.
Subsequently, general initial temperatures and specific heats were investigated -- see
\cite{BeGi94,vdB5,Ca45,McA93,McMe03a,PhJa90,Sav98a} and the references contained therein. Other boundary
conditions (Neumann, Zaremba, etc.) have been considered \cite{BeDeGi93,BGKK07}. The growth of the coefficients
$\beta_j^{\partial M}$ has also been of interest
\cite{x6,BGK11,TS11} -- see also \cite{x5} for related work on the heat trace asymptotics. The case where
$\psi_1$ is singular
$(\alpha_1>0)$ but $\psi_2$ is smooth
$(\alpha_2=0)$ was studied previously
\cite{vdB,vdB5}. The current paper is devoted to the study of the invariants $\beta_j^{\partial M}$ in the doubly singular case.

The special case of a ball of radius $a$ in
$\mathbb{R}^3$ is well understood. The following result was proved
in \cite{vdB9}.
\begin{theorem}
Let $B_a=\{x\in \R^3:|x|\le a\}$, and let $D$
be the Dirichlet Laplacian acting in $L^2(B_a)$.
If $\alpha_1<2, \alpha_2<2, \alpha_1+\alpha_2>3, J\in \N$ then there exist
coefficients $b_0, b_1,\cdots$ depending on $\alpha_1,\alpha_2$
only such that for $t \downarrow 0$
\begin{equation}\label{eqn-1}
\begin{array}{ll}
&Q(\delta^{-\alpha_1},\delta^{-\alpha_2},D)(t)\gronk \\
&=4\pi c_{\alpha_1, \alpha_2}a^2t^{(1-\alpha_1
-\alpha_2)/2}-4\pi(c_{\alpha_1-1,\alpha_2}+c_{\alpha_1,\alpha_2-1})at^{(2-\alpha_1 -\alpha_2)/2}\gronk\\
&\ \ +4\pi c_{\alpha_1-1,\alpha_2-1}t^{(3-\alpha_1
-\alpha_2)/2}+\displaystyle\sum_{j=0}^Jb_ja^{3-j-\alpha_1-\alpha_2}t^{j/2}+O(t^{(J+1)/2}),
\gronk\end{array}\end{equation} where
\begin{equation}\label{eqn-2}
\begin{array}{ll}
c_{\alpha_1, \alpha_2}=&2^{-\alpha_1 -\alpha_2}\pi^{-1/2}\Gamma((2-\alpha_1
-\alpha_2)/2)\gronk \\ & \times\displaystyle\int_0^1
(\rho^{-\alpha_1}+\rho^{-\alpha_2})((1-\rho)^{\alpha_1+\alpha_2-2}
  -(1+\rho)^{\alpha_1+\alpha_2-2})d\rho,\gronk
\end{array}\end{equation}
and
$$\begin{array}{ll}
b_0=-8\pi((\alpha_1+\alpha_2-1)(\alpha_1+\alpha_2-2)(\alpha_1+\alpha_2-3))^{-1},&b_1=0, \gronk \\
b_2=8\pi \alpha_1 \alpha_2
 ((\alpha_1+\alpha_2+1)(\alpha_1+\alpha_2)(\alpha_1+\alpha_2-1))^{-1}, &b_3=0.
\end{array}$$
\end{theorem}

It is convenient to use a standard formalism to describe the invariants $\beta_j^{\partial M}$ in the general setting. Let $D$ be an
operator of Laplace type acting on the space of smooth sections to some vector bundle $V$ over a Riemannian manifold $(M,g)$.
Choose a local system of coordinates
$(x^1,...,x^m)$ for $M$ and a local frame for $V$. We adopt the {\it Einstein
convention} and sum over repeated indices. Let $ds^2=g_{\mu\nu}dx^\mu\circ dx^\nu$ define the Riemannian metric and let
$g^{\mu\nu}$ be the inverse matrix where $1\le\mu,\nu\le m$. We may then express:
$$D=-(g^{\mu\nu}\operatorname{Id}\partial_{x_\mu}\partial_{x_\nu}+A^\mu\partial_{\mu}+B)$$
for suitably chosen endomorphisms $A^\mu$ and $B$ of $V$.
If $\nabla$
is a connection on
$V$, we use $\nabla$ and the Levi--Civita connection to covariantly differentiate tensors of all
types and let `;' denote multiple covariant differentiation. If $\psi_1$ is a section to $V$ which is smooth on
$\operatorname{int}(M)$, let
$\psi_{1;\mu\nu}$ be the components of
$\nabla^2\psi_1$. If $E$ is an auxiliary endomorphism of $V$, we define the associated {\it modified
Bochner Laplacian} by setting:
$$D(g,\nabla,E)\psi_1:=-g^{\mu\nu}\psi_{1;\nu\mu}-E\psi_1\,.$$
Let $\Gamma_{\mu\nu\sigma}$ and $\Gamma_{\mu\nu}{}^\sigma$ be the Christoffel symbols of the Levi-Civita connection.
Then (see, for example, the discussion in \cite{BeGi94}):
\begin{lemma}
If $D$ is an operator of Laplace type, then there exists a unique
connection $\nabla$ on $V$ and a unique endomorphism $E$ of $V$ so
that $D=D(g,\nabla,E)$. The connection $1$-form $\omega$ of
$\nabla$ and the endomorphism $E$ are given by
\begin{eqnarray*}
&&\omega_\mu=\textstyle{\textstyle\frac12}(g_{\mu\nu}A^\nu
+g^{\sigma\varepsilon}\Gamma_{\sigma\varepsilon\mu}
\operatorname{Id}),\\
&&E=B-g^{\mu\nu}(\partial_{x_\nu}\omega_{\mu}+\omega_{\mu}\omega_{\nu}
    -\omega_{\sigma}\Gamma_{\mu\nu}{}^\sigma)\,.
\end{eqnarray*}
\end{lemma}

The specific heat $\psi_2$ is a section to the dual vector bundle $\tilde V$. We use the dual connection
on $\tilde V$
to covariantly differentiate
$\psi_2$. Note that the connection
$1$ form
$\tilde\omega_\nu$ for $\tilde\nabla$ is the dual of $-\omega_\nu$. Thus
$$
  \tilde\nabla_{\partial_{x_\mu}}=\partial_{x_\mu}-
\textstyle{\textstyle\frac12}(g_{\mu\nu}\tilde A^\nu
+g^{\sigma\varepsilon}\Gamma_{\sigma\varepsilon\mu}\operatorname{Id})\,.
$$

Near the boundary, choose an orthonormal frame $\{e_1,...,e_m\}$
for the tangent bundle of $M$ so that $e_m$ is the inward unit
geodesic normal. Let indices $a,b$ range from $1$ to $m-1$ and
index the induced orthonormal frame $\{e_1,...,e_{m-1}\}$ for the
tangent bundle of the boundary. We let `$:$' denote the components
of tangential covariant  differentiation defined by $\nabla$ and
the Levi-Civita connection of the boundary. Let
$L_{ab}:=g(\nabla_{e_a}e_b,e_m)=\Gamma_{abm}$ be the components of
the second fundamental form. The difference between `$;$' and
`$:$' is then measured by $L$. Let $\tilde D$ be the dual operator
of Laplace type on $\tilde V$. The following relations will be
useful subsequently:
\begin{eqnarray*}
&&D\psi_1=-(\psi_{1:aa}+\psi_{1;mm}-L_{aa}\psi_{1;m}+E\psi_1),\\
&&\tilde D\psi_2=-(\psi_{2:aa}+\psi_{2;mm}-L_{aa}\psi_{2;m}+\tilde E\psi_2)\,.
\end{eqnarray*}

We expand $\psi_1$ and $\psi_2$ near the boundary of $M$ in the form:
$$\psi_1(y,\delta)\sim \delta^{-\alpha_1}\sum_{j=0}^\infty\psi_1^j\delta^j\ , \quad
  \psi_2(y,\delta)\sim \delta^{-\alpha_2}\sum_{j=0}^\infty\psi_2^j\delta^j,$$
where $\nabla_{e_m}\psi_1^j=0$ and $\tilde\nabla_{e_m}\psi_2^j=0$. We shall
usually be working with scalar operators and can choose local sections $s$ and $\tilde s$ so that $\nabla_{e_m}s=0$ and
$\tilde\nabla_{e_m}\tilde s=0$. We may then express $\psi_1=\Psi_1 s$, $\psi_1^i=\Psi_1^is$,
$\psi_2=\Psi_2 \tilde s$, and $\psi_2^i=\Psi_2^i\tilde s$ where $\Psi_1^j$ and $\Psi_2^j$ are smooth functions defined
on the boundary so that we have the modified Taylor series
$$\delta^{\alpha_i}\Psi_i(y,\delta)\sim\sum_{j=0}^\infty\Psi_i^j(y)\delta^j\quad\text{for}\quad i=1,2\,.$$
For the Laplacian, the bundles and connections under consideration
are trivial so this formalism is unnecessary. However, for more
general operators, the connections in question are not flat and
this formalism is essential. We shall be using the method of
``universal examples'' in what follows. It is a peculiar feature
of this method that even if we were only interested in the scalar
Laplacian for a smooth bounded domain in $\mathbb{R}^m$, it would
be necessary to deal with quite general operators as we shall see
presently while proving Lemma~\ref{lem-1.7} in
Section~\ref{sect-3}.

Let
$\operatorname{Ric}$ be the Ricci tensor of $M$, let $\tau$ be the scalar curvature of $M$, and let $dy$
be the Riemannian measure of $\partial M$.  Section~\ref{sect-3} is devoted to the proof of the following
result:
\begin{theorem}\label{thm-3}
Let $\alpha_1+\alpha_2\notin \Z$, $\alpha_1<2$, $\alpha_2<2$. Assume that Conjecture~\ref{conj-1.1}
holds. Let $c_{\alpha_1, \alpha_2}$ be as given in Equation~(\ref{eqn-2}). Then
\begin{eqnarray*}
&&\beta_0^{\partial M}=\textstyle\int_{\partial M}c_{\alpha_1,\alpha_2}\psi_1^0\psi_2^0dy,\\
&&\beta_1^{\partial M}=\textstyle\int_{\partial M}
\{c_{\alpha_1-1,\alpha_2}\psi_1^1\psi_2^0-
{\textstyle\frac12}\{c_{\alpha_1-1,\alpha_2}+c_{\alpha_1,\alpha_2-1}\}\psi_1^0\psi_2^0L_{aa}\\
&&\qquad+c_{\alpha_1,\alpha_2-1}\psi_1^0\psi_2^1\}dy,\\
&&\beta_2^{\partial M}=\textstyle\int_{\partial M}
\{c_{\alpha_1-2,\alpha_2}\psi_1^2\psi_2^0
-{\textstyle\frac12}(c_{\alpha_1-2,\alpha_2}+c_{\alpha_1-1,\alpha_2-1})L_{aa}\psi_1^1\psi_2^0\\
&&\qquad+
c_{\alpha_1,\alpha_2}E\psi_1^0\psi_2^0
+c_{\alpha_1,\alpha_2-2}\psi_1^0\psi_2^2
-{\textstyle\frac12}(c_{\alpha_1-1,\alpha_2-1}+c_{\alpha_1,\alpha_2-2})L_{aa}\psi_1^0\psi_2^1\\
&&\qquad+(-{\textstyle\frac14}c_{\alpha_1-2,\alpha_2}-{\textstyle\frac14}c_{\alpha_1,\alpha_2-2}
+{\textstyle\frac12}c_{\alpha_1,\alpha_2})
(L_{aa}L_{bb}+\hbox{\rm Ric}_{mm})\psi_1^0\psi_2^0\\
&&\qquad-c_{\alpha_1,\alpha_2}\psi_1^0{}_{:a}\rho_2^0{}_{0:a}
   +0\tau\psi_1^0\psi_2^0+c_{\alpha_1-1,\alpha_2-1}\psi_1^1\psi_2^1\\
&&\qquad+(\textstyle\frac18c_{\alpha_1-2,\alpha_2}+\frac18c_{\alpha_1,\alpha_2-2}
+\frac14c_{\alpha_1-1,\alpha_2-1}-\frac14c_{\alpha_1,\alpha_2})L_{aa}L_{bb}\psi_1^0\psi_2^0\}dy\,.
\end{eqnarray*}
\end{theorem}

For the ball $B_a$ in $\R^3$, $L_{aa}=2a^{-1}$, $L_{aa}L_{bb}=4a^{-2}$ and
$L_{aa}L_{bb}=2a^{-2}$. Hence
$$
\beta_1^{\partial B_a}=-4\pi a(c_{\alpha_1-1,\alpha_2}+c_{\alpha_1,\alpha_2-1}),
$$
and the coefficient of $t^{(2-\alpha_1-\alpha_2)/2}$ agrees with Equation~(\ref{eqn-1}) in Theorem~\ref{thm-3}.
Similarly, the next term in the series given by $\beta_2^{\partial B_a}$ is consistent with Theorem~\ref{thm-3}.

\medbreak We observe that the $\Gamma$- function in the expression for
$c_{\alpha_1, \alpha_2}$ which is given in Equation~(\ref{eqn-2}) has simple poles for
$\alpha_1+\alpha_2\in \{2,4,6,\cdots\}$. Furthermore the
integrand with respect to $\rho$ equals $0$ for
$\alpha_1+\alpha_2=2$. It is easily seen that this singularity is
removable. On the other hand the integral with respect to $\rho$
is finite only for $\alpha_1<2,\alpha_2<2$ and
$\alpha_1+\alpha_2>1$. This suggests that the $j^{\rm{th}}$ term
($j=1,2,3$) in Equation~(\ref{eqn-1}) will take a different form for
$\alpha_1+\alpha_2=j$. This is indeed the case for an interval in $\R$.
Let $a>0$, and let $\chi_1,\chi_2$ be non-negative $C^{\infty}$
functions on $\R^+$ defined by
$$
\chi_{1,2}(x)=\left\{\begin{array}{lll}1&\hbox{if}&0 \le x \le \epsilon_{1,2}, \\
0&\hbox{if}&\ x \ge \epsilon_{3,4},\\
\end{array}\right.
$$
where $0<\epsilon_1<\epsilon_3<a/2$ and
$0<\epsilon_2<\epsilon_4<a/2$. We shall establish the following result in Section~\ref{sect-2}:
\begin{theorem}\label{thm-4}
Let $\alpha_1<2,\alpha_2<2,
\alpha_1+\alpha_2=1$. If $t\downarrow 0$, then:
\begin{equation}\label{eqn-3}
\begin{array}{ll}
&\iint_{[0,a]^2}p_{[0,a]}(x_1,x_2;t)\chi_1(\delta(x_1))
\chi_2(\delta(x_2))\delta(x_1)^{-\alpha_1}\delta(x_2)^{-\alpha_2}dx_1dx_2\gronk
\\ &=\log(\epsilon^2/t)+\gamma+4\log(2^{1/2}-1)+4\log 2\gronk \\ &
\ +2\displaystyle\int_{[\epsilon,a/2]}\chi_1(x)\chi_2(x)x^{-1}dx
+\int_{[0,1]}dqq^{-1}(1+q^2)\gronk
\\ & \ \times\bigg\{((1+q)/(1-q))^{\alpha-1}+((1-q)/(1+q))^{\alpha}\gronk\\
&\ -2(1-q)(1+q^2)^{-1/2}\bigg\}+O(t^{1/2}\log t),\vphantom{\vrule height 15pt}
\end{array}\end{equation}
where $p_{[0,a]}(x_1,x_2;t), x_1\in [0,a],x_2\in[0,a], t>0$ is the
Dirichlet heat kernel for the interval $[0,a]$, $\gamma$ is
Euler's constant, and $\epsilon=\min\{\epsilon_1,\epsilon_2\}.$
\end{theorem}

We note that the $\epsilon$-dependence in the right hand side of
Equation~(\ref{eqn-3}) is fictitious. Since
$\chi_1(x)=\chi_2(x)=1$ for $0<x\le \epsilon,$ we have that
$$
\log(\epsilon^2/t)+2\int_{[\epsilon,a/2]}\chi_1(x)\chi_2(x)x^{-1}dx=2\int_{[\sqrt
t,a/2]}\chi_1(x)\chi_2(x)x^{-1}dx,
$$
which independent of $\epsilon$ for $0<t<\epsilon^2$. We also note
that the leading term in Theorem~\ref{thm-4} jibes with Theorem
1.4 (2) in \cite{vdB5} since the volume of the volume of the
boundary of the interval $[0,a]$ is equal to $2$. This supports
the following.
\begin{conjecture}\label{conj-2}
Let $M$ be a compact Riemannian manifold with smooth boundary $\partial M$, and let $\delta$ denote the
distance to the boundary. Let $\alpha_1<2,\alpha_2<2, \alpha_1+\alpha_2=1$, and
let $\chi_1$ and $\chi_2$ be smooth functions on $\R_+$ with
support contained in an interval $[0,b]$, and equal to $1$ in a
neighbourhood of $0$, and where $b$ is such that $\delta$ is
smooth on the collar $\partial M\times[0,b]$. If $t\downarrow 0$
then
$$
Q(\delta^{-\alpha_1}\chi_1\circ\delta,\delta^{-\alpha_2}\chi_2\circ\delta,D)(t)=2^{-1}\int_{\partial
M}dy\log t +o(\log t).
$$
\end{conjecture}

We note that Theorem~\ref{thm-4} and Conjecture~\ref{conj-2}
include the cases where either $1<\alpha_1<2$ or $1<\alpha_2<2$.
This requires more care in the proof of Theorem~\ref{thm-4} than
the case where both $\alpha_1<1$ and $\alpha_2<1$.
\section{The proof of Theorem~\ref{thm-4}}\label{sect-2}
The first step in the proof of Theorem~\ref{thm-4} is to reduce
the calculation on the interval $[0,a]$ to a calculation on the
half-line $\R^+=[0,\infty)$. We have the following:
\begin{lemma}
Let $\alpha_1<2,\alpha_2<2,
\alpha_1+\alpha_2=1$. If $t\downarrow 0$ then
\begin{equation}\label{eqn-4}
\begin{array}{ll}
&\iint_{[0,a]^2}p_{[0,a]}(x_1,x_2;t)\chi_1(\delta(x_1))\chi_2(\delta(x_2))
\delta(x_1)^{-\alpha_1}\delta(x_2)^{-\alpha_2}dx_1dx_2\gronk
\\
&=2\iint_{\R_+^2}p_{\R^+}(x_1,x_2;t)\chi_1(x_1)\chi_2(x_2)x_1^{-\alpha_1}
x_2^{-\alpha_2}dx_1dx_2\\
&\qquad+O(e^{-(1-\eta)\kappa^2/(4t)}),
\gronk\end{array}\end{equation}
where $\kappa=a-\epsilon_3-\epsilon_4$ and
$\eta=\max\{\alpha_1/2,\alpha_2/2\}$.
\end{lemma}
\proofx Without loss of generality we may assume that
$\alpha_1\ge \alpha_2$. We partition the region of integration
$[0,a]^2=\cup_{i=1}^5A_i$, where
$$\begin{array}{ll}
A_1=[0,\epsilon_3]\times[0,\epsilon_4],&
A_2=[0,\epsilon_3]\times[a-\epsilon_4,a],\gronk \\
A_3=[a-\epsilon_3,a]\times[0,\epsilon_4],&
A_4=[a-\epsilon_3,a]\times[a-\epsilon_4,a],\gronk \\
A_5=A\setminus(\cup_{i=1}^4A_i).\gronk\end{array}$$
The integrand in the left hand side of Equation~(\ref{eqn-4}) is identically
equal to $0$ on $A_5$, and this set does not contribute to the
integral. Since 
\begin{eqnarray*}
&&p_{[0,a]}(x_1,x_2;t)=p_{[0,a]}(a-x_1,a-x_2;t),\quad\text{and}\\
&&p_{[0,a]}(x_1,a-x_2;t)=p_{[0,a]}(a-x_1,x_2;t),
\end{eqnarray*}
the contributions  of $A_1$ and $A_2$ to the integral in the left hand
side of \eqref{eqn-4} are equal to the contributions of $A_4$ and $A_3$
respectively.
Since $|x_1-x_2|\ge \kappa$ for $(x_1,x_2)\in A_2$, we have
by monotonicity of the Dirichlet heat kernel that
\begin{equation}\label{eqn-5}
\begin{array}{ll}
p_{[0,a]}(x_1,x_2;t)&\le p_{\R^+}(x_1,x_2;t)\gronk \\ & =(4\pi
t)^{-1/2}\left(e^{-(x_1-x_2)^2/(4t)}-e^{-(x_1+x_2)^2/(4t)}\right)\gronk
\\ &=(4\pi
t)^{-1/2}e^{-(x_1-x_2)^2/(4t)}\left(1-e^{-x_1x_2/t}\right)\gronk \\ &
\le t^{-3/2}x_1x_2e^{-\kappa^2/(4t)}.
\gronk\end{array}\end{equation}
Hence the contribution from $A_2$ to the integral in the left hand
side of \eqref{eqn-4} is bounded from above by
$$t^{-3/2}e^{-\kappa^2/(4t)}\iint_{A_2}\chi_1(x_1)\chi_2(x_2)x_1^{1-\alpha_1}
x_2^{\alpha_1}dx_1dx_2=O(e^{-\eta\kappa^2/(4t)}).$$
The contribution from $A_1$ to the integral in the left hand side
of \eqref{eqn-4} is bounded from above by
\begin{eqnarray*}
&\iint_{A_1}&p_{\R^+}(x_1,x_2;t)\chi_1(x_1)\chi_2(x_2)x_1^{-\alpha_1}x_2^{-\alpha_2}dx_1dx_2
\\
&&=\iint_{\R_+^2}p_{\R^+}(x_1,x_2;t)\chi_1(x_1)\chi_2(x_2)x_1^{-\alpha_1}x
_2^{-\alpha_2}dx_1dx_2.
\end{eqnarray*}
This completes the proof of the upper bound.

To establish the lower bound we note that
\begin{eqnarray*}
&&\iint_{[0,a]^2}p_{[0,a]}(x_1,x_2;t)\chi_1(\delta(x_1))\chi_2(\delta(x_2))
\delta(x_1)^{-\alpha_1}\delta(x_2)^{-\alpha_2}dx_1dx_2
\\ &\ge&
2\iint_{A_1}p_{[0,a]}(x_1,x_2;t)\chi_1(x_1)\chi_2(x_2)x_1^{-\alpha_1}x_2^{-\alpha_2}dx_1dx_2.
\end{eqnarray*}
It is well known that the Dirichlet heat kernel for an open set
$\Omega \in \R^m$ has the following probabilistic representation.
$$
p_{\Omega}(x_1,x_2;t)=p_{\R^m}(x_1,x_2;t)\textup{Prob}_{x_1,x_2}[B(s)\in
\Omega,0<s<t],
$$
where $(B(s),0\le s \le t)$ is a Brownian bridge on $\R^m$. For
$\Omega=[0,a]\in \R$ and for $x\in [0,a],y\in [0,a]$ we have that
\begin{equation}\label{eqn-6}
\begin{array}{ll}
&p_{[0,a]}(x_1,x_2;t)=p_{\R}(x_1,x_2;t)\textup{Prob}_{x_1,x_2}[0<B(s)<a,0<s<t]\gronk
\\ &=p_{\R}(x_1,x_2;t)(\textup{Prob}_{x_1,x_2}[0<B(s),0<s<t]\gronk \\ &
\quad-\textup{Prob}_{x_1,x_2}[(0<B(s),0<s<t)\wedge
(\max_{0\le s\le t}B(s)\ge a)]).
\gronk\end{array}\end{equation}
By H\"older's inequality we have for $\eta \in(0,1)$
\begin{equation}\label{eqn-7}\begin{array}{ll}
&p_{\R}(x_1,x_2;t)\textup{Prob}_{x_1,x_2}[(0<B(s),0<s<t)\wedge (\max_{0\le
s\le t}B(s)\ge a)]\gronk \\ & \le
p_{\R}(x_1,x_2;t)(\textup{Prob}_{x_1,x_2}[0<B(s),0<s<t])^{\eta}\gronk
\\ & \hspace
{30mm}\times(1-\textup{Prob}_{x_1,x_2}[B(s)<a,0<s<t])^{1-\eta}\gronk
\\ &=(p_{\R^+}(x_1,x_2;t))^{\eta}(p_{\R}(x_1,x_2;t)-p_{(-\infty,a]}(x_1,x_2;t))^{1-\eta}.
\gronk\end{array}\end{equation}
Since
\begin{equation}\label{eqn-8}
p_{(-\infty,a]}(x_1,x_2;t)=(4\pi
t)^{-1/2}\left(e^{-(x_1-x_2)^2/(4t)}-e^{-(2a-x_1-x_2)^2/(4t)}\right),
\end{equation}
we have by Equation~(\ref{eqn-6}), Equation~(\ref{eqn-7}), and Equation~(\ref{eqn-8}) that
\begin{eqnarray*}
p_{[0,a]}(x_1,x_2;t)&\ge &p_{\R^+}(x_1,x_2;t)\\
&&\quad-(p_{\R^+}(x_1,x_2;t))^{\eta}(4\pi
t)^{-(1-\eta)/2}e^{-(1-\eta)(2a-x_1-x_2)^2/(4t)}.
\end{eqnarray*}
By the last inequality in Equation~(\ref{eqn-5})
\begin{eqnarray*}
&&(p_{\R^+}(x_1,x_2;t))^{\eta}(4\pi
t)^{-(1-\eta)/2}e^{-(1-\eta)(2a-x_1-x_2)^2/(4t)}\\
&&\qquad\qquad\le
t^{-\frac12-\eta}(x_1x_2)^{\eta}e^{-(1-\eta)(2a-x_1-x_2)^2/(4t)}.
\end{eqnarray*}
Integrating the above right hand side with respect to
$\chi_1(x_1)\chi_2(x_2)x_1^{-\alpha_1}x_2^{-\alpha_2}dx_1dx_2$ yields a bound
\begin{eqnarray*}
&&t^{-\frac12-\eta}\int_{[0,\epsilon_3]}\chi_1(x_1)x_1^{\eta-\alpha_1}dx_1
\int_{[0,\epsilon_4]}\chi_1(x_2)x_2^{\eta-\alpha_2}dx_2e^{-(1-\eta)(2a-x_1-x_2)^2/(4t)}\gronk
\\ &&\le t^{-\frac12-\eta}e^{-(1-\eta)(2a-\epsilon_3-\epsilon_4)^2/(4t)}\int_{[0,\epsilon_3]}
\chi_1(x_1)x_1^{\eta-\alpha_1}dx_1
\int_{[0,\epsilon_4]}\chi_1(x_2)x_2^{\eta-\alpha_2}dx_2\gronk \\ &&
=O(e^{-(1-\eta)(a-\epsilon_3-\epsilon_4)^2/(4t)}).
\end{eqnarray*}
Note that since $2>\alpha_1$ and $\eta=\alpha_1/2,$
$x_1^{\eta-\alpha_1}=x_1^{-\alpha_1/2}$ is integrable at $0$. Since
$1=\alpha_1+\alpha_2\le 2\alpha_1$ we have that $\alpha_1\ge
1/2>0$. Hence $x_2^{\eta-\alpha_2}=x_2^{-1+(3\alpha_1/2)}$ is also
integrable at $0$. This completes the proof of the lower bound.
\hfill\qedbox

In order to prove Theorem~\ref{thm-4} it clearly suffices to prove
the following.
\begin{lemma}
Let $\alpha_2\le \alpha_1<2,
\alpha_1+\alpha_2=1$. If $t\downarrow 0$ then
\begin{eqnarray}
&&\iint_{\R_+^2}p_{\R^+}(x_1,x_2;t)\chi_1(x_1)\chi_2(x_2)x_1^{-\alpha_1}x_2^{-\alpha_2}dx_1dx_2\nonumber\\
&=&2^{-1}\log(\epsilon^2/t)+2^{-1}\gamma+2\log(2^{1/2}-1)+2\log 2\nonumber\\
&&\quad+\int_{[\epsilon,a/2]}\chi_1(x)\chi_2(x)x^{-1}dx
+2^{-1}\int_{[0,1]}dqq^{-1}(1+q^2)\label{eqn-9}
\\
&&\qquad\times\bigg\{((1+q)/(1-q))^{\alpha-1}+((1-q)/(1+q))^{\alpha}\nonumber\\
&&\qquad-2(1-q)(1+q^2)^{-1/2}\bigg)+O(t^{1/2}\log t).\nonumber
\end{eqnarray}
\end{lemma}
\proofx
Define
\begin{eqnarray*}
&&C=\{(x_1,x_2)\in \R_+^2: x_1^2+x_2^2\ge \epsilon ^2,0\le x_1\le
\epsilon_3,0\le x_2 \le \epsilon_4\},\\
&&C_1=\{(x_1,x_2)\in C: |x_1-x_2|\le \sigma\},
\end{eqnarray*}
where $\sigma \in (0,\epsilon/5)$ will be chosen later on. The
left hand side of Equation~(\ref{eqn-9}) can be written as $B_1+B_2$, where
\begin{eqnarray}\label{eqn-10}
&&B_1=\iint_{\R_+^2\cap\{x_1^2+x_2^2<\epsilon^2\}}p_{\R^+}(x_1,x_2;t)x_1^{-\alpha_1}x_2^{-\alpha_2}dx_1dx_2,\\
&&B_2=\iint_{C}p_{\R^+}(x_1,x_2;t)\chi_1(x_1)\chi_2(x_2)x_1^{-\alpha_1}x_2^{-\alpha_2}dx_1dx_2.\nonumber
\end{eqnarray}

To estimate $B_2$ we first consider the contribution from the set
$C\setminus C_1$. We have by Equation~(\ref{eqn-5}) that
$$
p_{\R^+}(x_1,x_2;t)\le t^{-3/2}x_1x_2e^{-(x_1-x_2)^2/(4t)}\le
t^{-3/2}x_1x_2e^{-\sigma^2/(4t)}, (x_1,x_2)\in C\setminus C_1.
$$
Consequently,
\begin{equation}\label{eqn-11}
\iint_{C\setminus
C_1}p_{\R^+}(x_1,x_2;t)\chi_1(x_1)\chi_2(x_2)x_1^{-\alpha_1}x_2^{-\alpha_2}dx_1dx_2
\le Kt^{-3/2}e^{-\sigma^2/(4t)}
\end{equation}
 where
$$K=\iint_C\chi_1(x_1)\chi_2(x_2)x_1^{1-\alpha_1}x_2^{1-\alpha_2}dx_1dx_2.$$
On $C\cap \{|x_1-x_2|\le \epsilon/5\}$ we have that
$x_2\rightarrow \chi_2(x_2)x_2^{-\alpha_2}$ is $C^{\infty}$. Hence
there exists $L$ depending on $\epsilon,\alpha_2$ and on $\chi_2$
such that
$|\chi_2(x_2)x_2^{-\alpha_2}-\chi_2(_1)x_1^{-\alpha_2}|\le
L|x_1-x_2|$. It is easily seen that both $x_1\ge \epsilon/2$ and
$x_2\ge \epsilon/2$ on $C\cap \{|x_1-x_2|\le \epsilon/5\}$. Since
the Dirichlet heat kernel on $\R_+$ is bounded from above by
$t^{-1/2}$ we have that
\begin{equation}\label{eqn-12}
\begin{array}{ll}
\iint_{C_1}&p_{\R^+}(x_1,x_2;t)\chi_1(x_1)x_1^{-\alpha_1}L|x_1-x_2|dx_1dx_2\gronk
\\ & \le L(2/\epsilon)^{\alpha_1}t^{-1/2}\iint_{C_1}|x_1-x_2|dx_1dx_2\le
2aL(2/\epsilon)^{\alpha_1}t^{-1/2}\sigma^2.
\gronk\end{array}\end{equation} We now choose $\sigma^2$ as to
minimize $t^{-3/2}e^{-\sigma^2/(4t)}+t^{-1/2}\sigma^2$, i.e.
$$
\sigma^2=4t\log( t^{-2}).
$$
This gives that for $t$ sufficiently small the right hand sides of
Equation~(\ref{eqn-11}) and Equation~(\ref{eqn-12}) are $O(t^{1/2})$ and
$O(t^{1/2}\log(t^{-1}))$ respectively. We conclude that
\begin{equation}\label{eqn-13}
B_2=\iint_{C_1}p_{\R^+}(x_1,x_2;t)\chi_1(x_1)\chi_2(x_1)x_1^{-1}dx_1dx_2+O(t^{1/2}\log(t^{-1})).
\end{equation}
We now write
$$
C_1=(C_1\cap \{x_1^2\ge \epsilon^2/2\})\cup C_1 (\cap\{x_1^2< \epsilon^2/2\})=C_2 \cup C_3.
$$
Since $x_1\ge \epsilon/2$ on $C_1$ we have that the integrand in the first term in the right hand side of
Equation~(\ref{eqn-13}) is bounded by $2\epsilon^{-1}t^{-1/2}$. Hence
$$
\iint_{C_3}p_{\R^+}(x_1,x_2;t)\chi_1(x_1)\chi_2(x_1)x_1^{-1}dx_1dx_2\le 2\epsilon^{-1}t^{-1/2}|C_3|,
$$
where $|\cdot|$ denotes Lebesgue measure. It is easily seen that $|C_3|\le \sigma^2/2$. Consequently,
$$
0\le\iint_{C_3}p_{\R^+}(x_1,x_2;t)\chi_1(x_1)\chi_2(x_1)x_1^{-1}dx_1dx_2\le \epsilon^{-1}t^{-1/2}\sigma^2,
$$
and so the contribution from $C_3$ to the integral in Equation~(\ref{eqn-13}) is $O(t^{1/2}\log(t^{-1}))$.
Furthermore by monotonicity of the Dirichlet heat kernel 
$$p_{\R^+}(x_1,x_2;t)\le p_{\R}(x_1,x_2;t)\,.$$
Hence
\begin{eqnarray*}
&&\iint_{C_2}p_{\R^+}(x_1,x_2;t)\chi_1(x_1)\chi_2(x_1)x_1^{-1}dx_1dx_2\\
&\le& \iint_{\{x_1^2\ge
 \epsilon^2/2\}}p_{\R}(x_1,x_2;t)\chi_1(x_1)\chi_2(x_1)x_1^{-1}dx_1dx_2 \\
&=&\int_{[\epsilon/\sqrt2,a/2]}\chi_1(x_1)\chi_2(x_1)x_1^{-1}dx_1.
\end{eqnarray*}
To obtain a lower bound for the contribution from $C_2$ to the integral in Equation~(\ref{eqn-13})
we first observe that $(4\pi t)^{-1/2}e^{-(x_1+x_2)^2/(4t)}\le t^{-1/2}e^{-\epsilon^2/(4t)}$ and $x\ge
\epsilon/2$ for $(x_1,x_2)\in C_2$. Therefore
\begin{eqnarray*}
0&\le&\iint_{C_2}(4\pi t)^{-1/2}e^{-(x_1+x_2)^2/(4t)}\chi_1(x_1)\chi_2(x_1)x_1^{-1}dx_1dx_2\\
& \le&
2\epsilon^{-1}t^{-1/2}e^{-\epsilon^2/(4t)}|C_2|\le 2a^2\epsilon^{-1}t^{-1/2}e^{-\epsilon^2/(4t)}\\
&=&O(e^{-\epsilon^2/(5t)}).
\end{eqnarray*}
Finally
\begin{eqnarray*}
&&\iint_{C_2}p_{\R}(x_1,x_2;t)\chi_1(x_1)\chi_2(x_1)x_1^{-1}dx_1dx_2\gronk \\
&\ge&\iint_{\{x^2\ge \epsilon^2/2\}}p_{\R}(x_1,x_2;t)\chi_1(x_1)\chi_2(x_1)x_1^{-1}dx_1dx_2\\
&&\quad -\iint_{\{|x_1-x_2|\ge\sigma \}\cap \{x_1^2\ge \epsilon^2/2\}}p_{\R}(x_1,x_2;t)
\chi_1(x_1)\chi_2(x_1)x_1^{-1}dx_1dx_2\\
  &=&\int_{[\epsilon/\sqrt2,a/2]}\chi_1(x_1)\chi_2(x_1)x_1^{-1}dx_1\\
&&\quad -\iint_{\{|x_1-x_2|\ge\sigma \}\cap \{x_1^2\ge
\epsilon^2/2\}}p_{\R}(x_1,x_2;t)\chi_1(x_1)\chi_2(x_1)x_1^{-1}dx_1dx_2.
\end{eqnarray*}
Moreover
\begin{eqnarray*}
\int_{\{|x_1-x_2|\ge \sigma\}}p_{\R}(x_1,x_2;t)dx_2&\le& \int_{\{|x_1-x_2|\ge \sigma\}}(4\pi t)^{-1/2}
e^{-|x_1-x_2|\sigma/(4t)}dx_2
\nonumber \\
&=&4\pi^{-1/2}t^{1/2}\sigma^{-1}e^{-\sigma^2/(4t)}=O(t^2).
\end{eqnarray*}
Putting all this together gives that
\begin{eqnarray*}
B_2&=&\int_{[\epsilon/\sqrt2,a/2]}\chi_1(x)\chi_2(x)x^{-1}dx+O(t^{1/2}\log(t^{-1}))\\
&=&\int_{[\epsilon,a/2]}\chi_1(x)\chi_2(x)x^{-1}dx+2^{-1}\log
2+O(t^{1/2}\log(t^{-1})),
\end{eqnarray*}
since $\chi_1(x)\chi_2(x)x^{-1}=x^{-1}$ for $0<x\le \epsilon$.

In order to obtain the asymptotic behaviour of $B_1$ in
Equation~(\ref{eqn-10}), we introduce polar coordinates $x=(4t)^{1/2}\rho \cos
\theta, y=(4t)^{1/2}\rho \sin \theta$ to find that
\begin{eqnarray*}
B_1&=&\pi^{-1/2}\int_{[0,\pi/2]}d\theta(\cos \theta)^{-\alpha}(\sin
\theta)^{\alpha-1}\\
&&\quad\times\int_{[0,\epsilon/(4t)^{1/2}]}d\rho(e^{-\rho^2(1-\sin(2\theta))}-e^{-\rho^2(1+\sin(2\theta))}).
\end{eqnarray*}
A further change of variable $ \theta=\phi+\pi/4$ yields that
\begin{eqnarray*}
B_1&=&(2/\pi)^{1/2}\int_{[0,\pi/4]}d\phi\left(\frac{(\cos\phi+\sin
\phi)^{\alpha-1}}{(\cos\phi-\sin
\phi)^{\alpha}}+\frac{(\cos\phi-\sin
\phi)^{\alpha-1}}{(\cos\phi+\sin \phi)^{\alpha}}\right)\\
& & \ \times\int_{[0,\epsilon/(4t)^{1/2}]}d\rho\left(e^{-2\rho^2(\sin \phi)^2}-e^{-2\rho^2(\cos
\phi)^2}\right)=B_3+B_4+B_5,
\nonumber\end{eqnarray*}
where
\begin{eqnarray}
B_3&=&(2/\pi)^{1/2}\int_{[0,\pi/4]}d\phi\left(\frac{(\cos\phi+\sin
\phi)^{\alpha-1}}{(\cos\phi-\sin
\phi)^{\alpha}}+\frac{(\cos\phi-\sin
\phi)^{\alpha-1}}{(\cos\phi+\sin \phi)^{\alpha}}-2\right)\nonumber
\\&&\quad\times\int_{[0,\infty)}d\rho\left(e^{-2\rho^2(\sin \phi)^2}-e^{-2\rho^2(\cos
\phi)^2}\right)\nonumber\\
&=&2^{-1}\int_{[0,\pi/4]}d\phi(\cos\phi)^{-1}(\sin\phi)^{-1}\nonumber
\\ &&\quad  \times\left(\frac{(\cos\phi+\sin \phi)^{\alpha-1}}{(\cos\phi-\sin
\phi)^{\alpha-1}}+\frac{(\cos\phi-\sin
\phi)^{\alpha}}{(\cos\phi+\sin
\phi)^{\alpha}}-2(\cos\phi-\sin\phi)\right)\label{eqn-14}\\
&=&2^{-1}\int_{[0,1]}dqq^{-1}(1+q^2)
\times\bigg\{((1+q)/(1-q))^{\alpha-1}\nonumber
\\&&\qquad\qquad+((1-q)/(1+q))^{\alpha}-2(1-q)(1+q^2)^{-1/2}\bigg\},\nonumber
\end{eqnarray}
\begin{eqnarray*}
B_4&=&-(2/\pi)^{1/2}\int_{[0,\pi/4]}d\phi\\
&&\quad\times\left\{\frac{(\cos\phi+\sin
\phi)^{\alpha-1}}{(\cos\phi-\sin
\phi)^{\alpha}}+\frac{(\cos\phi-\sin
\phi)^{\alpha-1}}{(\cos\phi+\sin \phi)^{\alpha}}-2\right\}\\
&&\quad\times\int_{[\epsilon/(4t)^{1/2},\infty)}d\rho\left(e^{-2\rho^2(\sin \phi)^2}-e^{-2\rho^2(\cos
\phi)^2}\right),
\end{eqnarray*}
and
\begin{equation}\label{eqn-15}
B_5=(8/\pi)^{1/2}\int_{[0,\pi/4]}d\phi\int_{[0,\epsilon/(4t)^{1/2}]}d\rho\left(e^{-2\rho^2(\sin
\phi)^2}-e^{-2\rho^2(\cos \phi)^2}\right).
\end{equation}
We have used the standard change of variables $\tan \phi=q$ to
obtain the last identity in Equation~(\ref{eqn-14}).

In order to find the asymptotic behaviour of $B_5$ as $t\downarrow
0$ we first consider the contribution of the second term in the
integrand with respect to $\rho$ in \eqref{eqn-15}, and write
\begin{eqnarray*}
&&-(8/\pi)^{1/2}\int_{[0,\pi/4]}d\phi\int_{[0,\epsilon/(4t)^{1/2}]}d\rho
e^{-2\rho^2(\cos \phi)^2}\\
&=&-\int_{[0,\pi/4]}d\phi(\cos\phi)^{-1}
+(8/\pi)^{1/2}\int_{[0,\pi/4]}d\phi\int_{[\epsilon/(4t)^{1/2},\infty)}d\rho
e^{-2\rho^2(\cos \phi)^2}\\
&=&\log(2^{1/2}-1)+O(e^{-\epsilon^2/(5t)}).
\end{eqnarray*}
The contribution of the first term in the integrand with respect
to $\rho$ in \eqref{eqn-15} is calculated as follows:
\begin{equation}\label{eqn-16}
\begin{array}{ll}
&\ \ \
(8/\pi)^{1/2}\int_{[0,\pi/4]}d\phi\int_{[0,\epsilon/(4t)^{1/2}]}d\rho
e^{-2\rho^2(\sin \phi)^2}\gronk \\ &
=(8/\pi)^{1/2}\int_{[0,\pi/4]}d\phi\int_{[0,\epsilon/(4t)^{1/2}]}d\rho
e^{-2\rho^2\phi^2}\\& \ \ \ +\int_{[0,\pi/4]}d\phi((\sin\phi)^{-1}-\phi^{-1})\gronk
\\ & \ \ \
+(8/\pi)^{1/2}\int_{[0,\pi/4]}d\phi\int_{[\epsilon/(4t)^{1/2},\infty)}d\rho\left(e^{-2\rho^2\phi^2}
-e^{-2\rho^2(\sin\phi)^2}\right).
\gronk\end{array}\end{equation}
The third term in the right hand side of Equation~(\ref{eqn-16}) is
$O(e^{-\epsilon^2/(5t)})$. The second term in the right hand side
of \eqref{eqn-16} is equal to $\log(2^{1/2}-1)+3\log 2-\log \pi$. The
first term in the right hand side of \eqref{eqn-16} equals
\begin{eqnarray*}
&&
(4/\pi)^{1/2}\int_{[0,\pi\epsilon/(32t)^{1/2}]}d\phi\int_{[0,\phi]}d\rho
e^{-\rho^2}\\
&=&(4/\pi)^{1/2}\left(\log(\pi\epsilon/(32t)^{1/2}\right)\int_{[0,\pi\epsilon/(32t)^{1/2}]}d\rho
e^{-\rho^2}\\
&&-(4/\pi)^{1/2}\int_{[0,\pi\epsilon/(32t)^{1/2}]}d\phi(\log
\phi)e^{-\phi^2}\\
&=&2^{-1}\log(\epsilon^2/t)+\log\pi+2^{-1}\gamma-3\cdot2^{-1}\log
2+O(e^{-\epsilon^2/(5t)}),
\end{eqnarray*}
where we have used Equation~(4.333) in \cite{GR} together with
$$
\int_{[0,\pi\epsilon/(32t)^{1/2}]}d\phi(\log
\phi)e^{-\phi^2}=\int_{[0,\infty)}d\phi(\log
\phi)e^{-\phi^2}+O(e^{-\epsilon^2/(5t)}).
$$
We find that
$$
B_5=2^{-1}\log(\epsilon^2/t)+2^{-1}\gamma+2\log(2^{1/2}-1)+3\cdot2^{-1}\log
2+O(e^{-\epsilon^2/(5t)}).
$$

In order to estimate $B_4$ we first note that by expanding
$\sin\phi$ and $\cos\phi$ around $0$ we have that
$$
\frac{(\cos\phi+\sin \phi)^{\alpha-1}}{(\cos\phi-\sin
\phi)^{\alpha}}+\frac{(\cos\phi-\sin
\phi)^{\alpha-1}}{(\cos\phi+\sin \phi)^{\alpha}}-2=O(\phi^2).
$$
Furthermore for $\phi \in [0,\pi/4]$,
\begin{eqnarray*}
0&\le&
\int_{[\epsilon/(4t)^{1/2},\infty)}d\rho\left(e^{-2\rho^2(\sin
\phi)^2}-e^{-2\rho^2(\cos \phi)^2}\right)\gronk \\ &
\le&(4t)^{1/2}\epsilon^{-1}\int_{[\epsilon/(4t)^{1/2},\infty)}d\rho\left(e^{-2\rho^2(\sin
\phi)^2}-e^{-2\rho^2(\cos \phi)^2}\right)\gronk \\
&\le&(4t)^{1/2}\epsilon^{-1}\int_{[\epsilon/(4t)^{1/2},\infty)}d\rho\rho\left(e^{-2\rho^2(\sin
\phi)^2}-e^{-2\rho^2(\cos \phi)^2}\right)\gronk \\
&\le&(4t)^{1/2}\epsilon^{-1}\int_{[0,\infty)}d\rho\rho\left(e^{-2\rho^2(\sin
\phi)^2}-e^{-2\rho^2(\cos \phi)^2}\right)\gronk \\
&=&t^{1/2}\epsilon^{-1}\left((\sin\phi)^{-2}-(\cos\phi)^{-2}\right).
\end{eqnarray*}
Hence
\begin{eqnarray*}
&&|B_4|\le
\pi^{-1/2}t^{1/2}\epsilon^{-1}\int_{[0,\pi/4]}d\phi(\sin\phi)^{-2}(\cos\phi)^{-2}
\\\quad &&\times\left|\frac{(\cos\phi+\sin \phi)^{\alpha}}{(\cos\phi-\sin
\phi)^{\alpha-1}}+\frac{(\cos\phi-\sin
\phi)^{\alpha}}{(\cos\phi+\sin
\phi)^{\alpha-1}}-2((\cos\phi)^2-(\sin\phi)^2)\right|.\nonumber
\end{eqnarray*}
We see that the integral with respect to $\phi$ converges both at
$\phi=0$ and at $\phi=\pi/4$. We conclude that $B_4=O(t^{1/2})$.\hfill\qedbox\medbreak

\section{The proof of Theorem~\ref{thm-3}}\label{sect-3}
We shall assume that $\Re(\alpha_1)<<0$ and $\Re(\alpha_2)<<0$ and then apply analytic continuation to establish the general case.
We shall also assume that
$\alpha_1+\alpha_2\notin\mathbb{Z}$ to ensure that the interior and the boundary terms do not interact. The invariants
$\beta_j^{\partial M}$ are given by local formula. Standard arguments using dimensional analysis yield the following result;
as these arguments are by now standard (see, for example, the discussion in \cite{BeGi94}), we omit details in the interests of
brevity.

\begin{lemma}\label{lem-7}
There exist universal constants $\varepsilon_{\alpha_1,\alpha_2}^i$  so
that:
\medbreak$
\beta_0^{\partial M}=\int_{\partial
M}\varepsilon_{\alpha_1,\alpha_2}^0\langle\psi_1^0,\psi_2^0\rangle dy$,
\medbreak$
\beta_1^{\partial M}=\int_{\partial M}\left\{
\varepsilon_{\alpha_1,\alpha_2}^1\langle\psi_1^1,\psi_2^0\rangle+\varepsilon_{\alpha_1,\alpha_2}^2
\langle L_{aa}\psi_1^0,\psi_2^0\rangle
+\varepsilon^3_{\alpha_1,\alpha_2}\langle\psi_1^0,\psi_2^1\rangle\right\}dy$,
\medbreak$
\beta_2^{\partial M}
    =\int_{\partial M}\{\varepsilon_{\alpha_1,\alpha_2}^4\langle\psi_1^2,\psi_2^0\rangle
+\varepsilon_{\alpha_1,\alpha_2}^5\langle
L_{aa}\psi_1^1,\psi_2^0\rangle
    +\varepsilon_{\alpha_1,\alpha_2}^6\langle E\psi_1^0,\psi_2^0\rangle$
$+\varepsilon_{\alpha_1,\alpha_2}^7\langle\psi_1^0,\psi_2^2\rangle$
\medbreak$\quad
   +\varepsilon_{\alpha_1,\alpha_2}^8\langle L_{aa}\psi_1^0,\psi_2^1\rangle
    +\varepsilon_{\alpha_1,\alpha_2}^9\langle\hbox{\rm Ric}_{mm}\psi_1^0,\psi_2^0\rangle
+\varepsilon_{\alpha_1,\alpha_2}^{10}\langle L_{aa}L_{bb}\psi_1^0,\psi_2^0\rangle$
\medbreak$\quad
+\varepsilon_{\alpha_1,\alpha_2}^{11}\langle
L_{ab}L_{ab}\psi_1^0,\psi_2^0\rangle
+\varepsilon_{\alpha_1,\alpha_2}^{12}\langle \psi_1^0{}_{:a},\psi_{2:a}^0\rangle
+\varepsilon_{\alpha_1,\alpha_2}^{13}\langle\tau\psi_1^0,\psi_2^0\rangle$
\medbreak$\quad+
\varepsilon_{\alpha_1,\alpha_2}^{14}\langle\psi_1^1,\psi_2^1\rangle\}dy$.
\end{lemma}

\begin{remark}\rm We note that $\varepsilon_{\alpha_1,\alpha_2}^0=c_{\alpha_1,\alpha_2}$ is given by
Equation~(\ref{eqn-2}).\end{remark}

There is a basic symmetry which is useful. Let $e^{-tD}$ denote the fundamental solution of the Dirichlet Laplacian and let
$\tilde D$ be the dual operator on the dual vector bundle $\tilde V$. The lemma below follows immediately from the
identity
$$Q(\psi_1,\psi_2,D)(t)=
\int_M\langle e^{-tD}\psi_1,\psi_2\rangle dx
=\int_M\langle\psi_1,e^{-t\tilde D}\psi_2\rangle dx=Q(\psi_2,\psi_1,\tilde D)(t)\,.
$$
\begin{lemma} Adopt the notation of Lemma~\ref{lem-7}.
$$\begin{array}{llllll}
\varepsilon_{\alpha_1,\alpha_2}^0=\varepsilon_{\alpha_2,\alpha_1}^0,&
\varepsilon_{\alpha_1,\alpha_2}^1=\varepsilon_{\alpha_2,\alpha_1}^3,&
\varepsilon_{\alpha_1,\alpha_2}^2=\varepsilon_{\alpha_2,\alpha_1}^2,&
\varepsilon_{\alpha_1,\alpha_2}^4=\varepsilon_{\alpha_2,\alpha_1}^7,\\
\varepsilon_{\alpha_1,\alpha_2}^5=\varepsilon_{\alpha_2,\alpha_1}^8,&
\varepsilon_{\alpha_1,\alpha_2}^6=\varepsilon_{\alpha_2,\alpha_1}^6,&
\varepsilon_{\alpha_1,\alpha_2}^9=\varepsilon_{\alpha_2,\alpha_1}^9,&
\varepsilon_{\alpha_1,\alpha_2}^{10}=\varepsilon_{\alpha_2,\alpha_1}^{10},\\
\varepsilon_{\alpha_1,\alpha_2}^{11}=\varepsilon_{\alpha_2,\alpha_1}^{11},&
\varepsilon_{\alpha_1,\alpha_2}^{12}=\varepsilon_{\alpha_2,\alpha_1}^{12},&
\varepsilon_{\alpha_1,\alpha_2}^{13}=\varepsilon_{\alpha_2,\alpha_1}^{13},&
\varepsilon_{\alpha_1,\alpha_2}^{14}=\varepsilon_{\alpha_2,\alpha_1}^{14}.
\end{array}$$
\end{lemma}

Next, we consider some product formulae:
\begin{lemma}\label{lem-1.5}
 Suppose that $M=M_1\times M_2$, that $g_M=g_{M_1}+g_{M_2}$, that $\partial {M_1}=\emptyset$, and
that
$D_M=D_{M_1}+D_{M_2}$ where $D_{M_1}$ and $D_{M_2}$ are scalar operators of Laplace type on ${M_1}$  and on
${M_2}$, respectively. Suppose that
$\psi_1^{M}=\psi_1^{M_1}\psi_1^{M_2}$ and $\psi_2^{M}=\psi_2^{M_1}\psi_2^{M_2}$ decompose similarly. Then
\begin{enumerate}
\item[\textup{(a)}]
$\beta(\psi_1^{M},\psi_2^{M},D_M)(t)=\beta(\psi_1^{M_1},\psi_2^{M_1},D_{M_1})(t)\cdot\beta(\psi_1^{M_2},\psi_2^{M_2},D_{M_2})(t)$.
\smallbreak\item[\textup{(b)}] $\int_{\partial
M}\beta_{k,{\alpha_1,\alpha_2}}^{\partial
M}(\psi_1^{M},\psi_2^{M},D_M)dy=
\sum_{2n+j=k}\frac{(-1)^n}{n!}\int_{M_1}\langle\psi_1^{M_1},(\tilde
D_{M_1})^n\psi_2^{M_1}\rangle dx_{M_1}$
\smallbreak$\quad\qquad\qquad\qquad\qquad\qquad\qquad\times
\int_{\partial M_2}\beta_{j,{\alpha_1,\alpha_2}}^{\partial
M_2}(\psi_1^{M_2},\psi_2^{M_2},D_{M_2})dy_{M_2}$.
\smallbreak\item[\textup{(c)}] The universal constants
$\varepsilon_{\alpha_1,\alpha_2}^i$ are dimension free.
\smallbreak\item[\textup{(d)}]
$\varepsilon_{\alpha_1,\alpha_2}^6=\varepsilon_{\alpha_1,\alpha_2}^0$,
$\varepsilon_{\alpha_1,\alpha_2}^{13}=0$, and
$\varepsilon_{\alpha_1,\alpha_2}^{12}=-\varepsilon_{\alpha_1,\alpha_2}^0$.
\end{enumerate}\end{lemma}

\medbreak\noindent{\bf Proof.} Assertion (a) follows from the identity $e^{-tD_M}
=e^{-tD_{M_1}}e^{-tD_{{M_2}}}$ and Assertion (b) follows from
Assertion (a). If we take $M_1=S^1$, $D_{M_1}=-\partial_\theta^2$, $\psi_1^{M_1}=1$,
and $\psi_2^{M_1}=1$, we have that
$\beta(\psi_1^{M_1},\psi_2^{M_1},D_{M_1})(t)=2\pi$. This then yields the identity
$$\int_{\partial M}\beta_{k,{\alpha_1,\alpha_2}}^{\partial M}(\psi_1^{M_2},\psi_2^{M_2},D)dy
=2\pi\int_{\partial M_2}\beta_{k,{\alpha_1,\alpha_2}}^{\partial
M_2}(\psi_1^{M_2},\psi_2^{M_2},D_{M_2})dy_2\,.$$ Assertion (c) now follows.
We take $M_2=[0,1]$ and
$D_2=-\partial_r^2$. We take
\begin{eqnarray*}
&&\psi_1^{M_2}=\psi_2^{M_2}=0\quad\text{near}\quad r=1,\\
&&\psi_2^{M_2}=r^{-\alpha_2}\quad\text{and}\quad\psi_1^{M_2}=r^{-\alpha_1}\quad\text{near}\quad r=0\,.
\end{eqnarray*}
Since the structures on $M_2$ are flat, we have $\psi_1^k=\psi_2^k=0$ for $k>0$ while
$$\psi_2^0=\psi_1^0=\left\{\begin{array}{lll}0&\text{at}&r=1\\ 1&\text{at}&r=0\end{array}\right\}\,.$$
Consequently,
$$
\beta_k^{\partial M_2}(\psi_1^{M_2},\psi_2^{M_2},D_{M_2})(r)=
\left\{\begin{array}{lll}
0&\text{if}&r=1\quad\text{and}\quad k\ge0\\
0&\text{if}&r=0\quad\text{and}\quad k>0\\
\varepsilon_{\alpha_1,\alpha_2}^0&\text{if}&r=0\quad\text{and}\quad k=0
\end{array}\right\}.$$
As the second fundamental form vanishes, the distinction between
`;' and `:' disappears, and we have $\tilde
D_1\psi_2^{M_1}=-(\psi_{2;aa}^{M_1}+\tilde E\psi_2)$. Calculating
on the interior then implies that
$$
\beta_{2}(\psi_1^{M_1},\psi_2^{M_1},D_{M_1})=\int_{M_1}\langle\psi_1^{M_1},\psi^{M_1}_{2;aa}+\tilde E\psi_2^{M_1}\rangle dx_1\,.
$$
We may therefore use Assertion (b) to derive the following identity from which Assertion~(e) will follow:
\medbreak\hfill
$\displaystyle\int_{\partial M}\beta_{2,{\alpha_1,\alpha_2}}^{\partial M}(\psi_1^{M},\psi_2^{M},D_M)dy
=\varepsilon_{\alpha_1,\alpha_2}^0\int_{M_1}\langle\psi_1^{M_1},\psi^{M_1}_{2;aa}+\tilde
E\psi_2^{M_1}\rangle dx_1$.
\hfill\qedbox\medbreak

We continue our study by index shifting:
\begin{lemma}
$$\begin{array}{llll}
\varepsilon_{\alpha_1,\alpha_2}^1=\varepsilon_{\alpha_1-1,\alpha_2}^0,&
\varepsilon_{\alpha_1,\alpha_2}^4=\varepsilon_{\alpha_1-2,\alpha_2}^0,&
\varepsilon_{\alpha_1,\alpha_2}^5=\varepsilon_{\alpha_1-1,\alpha_2}^2,\\
\varepsilon_{\alpha_1,\alpha_2}^3=\varepsilon^0_{\alpha_1,\alpha_2-1},&
\varepsilon_{\alpha_1,\alpha_2}^7=\varepsilon^0_{\alpha_1,\alpha_2-2},&
\varepsilon_{\alpha_1,\alpha_2}^8=\varepsilon_{\alpha_1,\alpha_2-1}^2,\\
\varepsilon_{\alpha_1,\alpha_2}^{14}=\varepsilon_{\alpha_1-1,\alpha_2-1}^0.
\end{array}$$
\end{lemma}

\medbreak\noindent{\bf Proof.} We assume $\psi_1$ and $\psi_2$ have compact support near the boundary of
$M$.  We set
$\tilde\psi_1:=(\delta^{n_1}\psi)\delta^{-\alpha_1-n_1}$ and $\tilde\psi_2:=(\delta^{n_1}\psi_2)\delta^{-\alpha_2-n_2}$ for
$n_i\in\mathbb{N}$. We compute:
\begin{eqnarray*}
&&\sum_kt^{(1+k-n_1-\alpha_1-n_2-\alpha_2)/2}\int_{\partial
M}\beta_{k,n_1+\alpha_1,n_2+\alpha_2}(\tilde\psi_1,\tilde\psi_2,D)dy\\
&\sim&\sum_{\ell}t^{(1+\ell-\alpha_1-\alpha_2)/2}\int_{\partial
M}\beta_{\ell,\alpha_1,\alpha_2}(\psi_1,\psi_2,D)dy\,.
\end{eqnarray*}
We set $k=\ell+n_1+n_2$ and equate powers of $t$ to see
$$\beta_{\ell+n_1+n_2,n_1+\alpha_1,n_2+\alpha_2}(\tilde\psi_1,\tilde\psi_2,D)
=\beta_{\ell,\alpha_1,\alpha_2}(\psi_1,\psi_2,D)\,.$$
 Note that
$\tilde\psi_1^{\mu+n_1}=\psi_1^\mu$ and
$\tilde\psi_2^{\nu+n_2}=\psi_2^\nu$. The desired result now follows by taking $(n_1,n_2)=(1,0)$, $(0,1)$,
$(2,0)$, $(1,1)$, and $(0,2)$.\hfill\qedbox\medbreak

\begin{lemma}\label{lem-1.7}
Let $\mathbb{T}^{m-1}$ denote the torus with periodic parameters
$(y_1,...,y_{m-1})$ and let $M:=\mathbb{T}^{m-1}\times[0,1]$.
Let $f_a\in C^\infty([0,1])$ have compact support near $r=0$ with $f_a(0)=0$. Let  $\Theta(r)\in
C^\infty([0,1])$ have compact support near $r=0$ with $\Theta\equiv1$ near $r=0$. Let
$\delta_a\in\mathbb{R}$. Set
$$\begin{array}{ll}
ds^2_M=\textstyle\sum_ae^{2f_a(r)}dy_a\circ dy_a+dr\circ dr,&
\psi_2:=\Theta(r)e^{-\sum_af_a(r)}r^{-\alpha_2},\\
D_M:=-\textstyle\sum_ae^{-2f_a(r)}(\partial_{y_a}^2+\delta_a\partial_{y_a})-\partial_r^2,
&\psi_1:=\Theta(r)r^{-\alpha_1}\,.\vphantom{\vrule height 12pt}
\end{array}$$
\begin{enumerate}
\item[\textup{(a)}] If $k>0$, then {$\int_{\partial M}\beta_{k,{\alpha_1,\alpha_2}}^{\partial M}(\psi_1,\psi_2,D_M)dy=0$}.
\smallbreak\item[\textup{(b)}]
$-{\textstyle\frac12}\varepsilon_{\alpha_1,\alpha_2}^1-\varepsilon_{\alpha_1,\alpha_2}^2
-{{\textstyle\frac12}\varepsilon_{\alpha_1,\alpha_2}^3}=0$.
\smallbreak\item[\textup{(c)}]
$-{\textstyle\frac14}(\varepsilon_{\alpha_1,\alpha_2}^6
+\varepsilon_{\alpha_1,\alpha_2}^{12})=0$.
\smallbreak\item[\textup{(d)}]
$-{\textstyle\frac14}\varepsilon_{\alpha_1,\alpha_2}^4
+{\textstyle\frac12}\varepsilon_{\alpha_1,\alpha_2}^6-{\textstyle\frac14}\varepsilon_{\alpha_1,\alpha_2}^7
-\varepsilon_{\alpha_1,\alpha_2}^9=0$.
\smallbreak\item[\textup{(e)}]
${\textstyle\frac18}\varepsilon_{\alpha_1,\alpha_2}^4
+{\textstyle\frac12}\varepsilon_{\alpha_1,\alpha_2}^5+{\textstyle\frac14}\varepsilon_{\alpha_1,\alpha_2}^6
+{\textstyle\frac18}\varepsilon_{\alpha_1,\alpha_2}^7
+{\textstyle\frac12}\varepsilon_{\alpha_1,\alpha_2}^{8}
+\varepsilon_{\alpha_1,\alpha_2}^{10}{+\frac14\varepsilon_{\alpha_1,\alpha_2}^{14}}=0$.
\smallbreak\item[\textup{(f)}]
$\textstyle-\varepsilon_{\alpha_1,\alpha_2}^9+\varepsilon_{\alpha_1,\alpha_2}^{11}=0$.
\end{enumerate}
\end{lemma}

\medbreak\noindent{\bf Proof.} We use $-\partial_r^2$ on $[0,1]$
and $D_M$ on $M$. Since $\Theta$ vanishes near $r=1$, this
boundary component plays no role. Let $u(r;t)$ be the solution of
the heat equation on $[0,1]$ with Dirichlet boundary conditions
and initial temperature $\psi_1$. The parameter $r$ is the
geodesic distance to the boundary near $r=0$. Since the problem
decouples, $u(r;t)$ is also the solution of the heat equation on
$M$ with Dirichlet boundary conditions. The Riemannian measure
$$dx=\sqrt{\det g_{ij}}dydr=e^{\sum_af_a}dydr\,.$$
As $\psi_2=\Theta e^{-\sum_af_a}r^{-\alpha_2}$, $\psi_2 dx=\Theta r^{-\alpha_2}dydr$. Since
$\operatorname{vol}(\mathbb{T}^{m-1})=(2\pi)^{m-1}$,
\begin{eqnarray*}
&&Q(\psi_1,\psi_2,D)(t)=
\int u(r;t)\psi_2 dx=(2\pi)^{m-1}\int_0^1u(r;t)\Theta(r)r^{-\alpha_2}dr\\
&=&(2\pi)^{m-1}\beta(\Theta r^{-\alpha_1},\Theta
r^{-\alpha_2},-\partial_r^2)(t)\,.
\end{eqnarray*}
The structures are flat on $[0,1]$. Since $\Theta$ vanishes identically near $r=1$ and  $\Theta$ is
identically $1$ near $r=0$, only the term $\beta_0$ is relevant in computing the boundary terms; the
$\beta_{k,\alpha_1,\alpha_2}$ vanish for $k\ge1$.

To apply Assertion (a), we must determine the relevant tensors. We have:
$$
\begin{array}{ll}
\Gamma_{abm}=-f_a^\prime\delta_{ab}e^{2f_a},&\Gamma_{ab}{}^m=-f_a^\prime e^{2f_a}\delta_{ab},\\
\Gamma_{amb}=f_a^\prime\delta_{ab}e^{2f_a} ,&\Gamma_{am}{}^b=f_a^\prime\delta_{a,b},
   \vphantom{\vrule height 12pt}\\
L_{ab}=\Gamma_{ab}{}^m|_{\partial M}=-f_a^\prime\delta_{ab},&\vphantom{\vrule height 12pt}\\
\textstyle\omega_a={\textstyle\frac12}e^{2f_a}\delta_a,&\tilde\omega_a=-\omega_a
   =-{\textstyle\frac12}e^{2f_a}\delta_a,
   \vphantom{\vrule height 12pt}\\
\omega_m=-{\textstyle\frac12}\sum_af_a^\prime,&
\tilde\omega_m=-\omega_m={\textstyle\frac12}\sum_af_a^\prime\,.\vphantom{\vrule height 12pt}
\end{array}$$
Consequently:
\medbreak\qquad
$R_{ambm}=g((\nabla_a\nabla_m-\nabla_m\nabla_a)e_b,e_m)
=\Gamma_{ac}{}^m\Gamma_{mb}{}^c-\partial_m\Gamma_{ab}{}^m$
\smallbreak\qquad\quad
$=\{-(f_a^\prime)^2+f_a^{\prime\prime}+2(f_a^\prime)^2\}e^{2f_a}\delta_{ab}$,
\medbreak\qquad
$\operatorname{Ric}_{mm}=-\textstyle\sum_a\left\{f_a^{\prime\prime}+(f_a^\prime)^2\right\}$,
\medbreak\qquad
$E|_{\partial M}=-\partial_m\omega_m-\omega_a^2-\omega_m^2+\omega_m\Gamma_{aa}{}^m$
\smallbreak\qquad\quad
$=\textstyle\frac12\sum_af_a^{\prime\prime}-\frac14\sum_a\delta_a^2-\frac14\sum_{a,b}f_a^\prime
f_b^\prime+\frac12\sum_{a,b}f_a^\prime f_b^\prime$
\medbreak\qquad\quad
$=\textstyle\frac12\sum_af_a^{\prime\prime}-\frac14\sum_a\delta_a^2+\frac14\sum_{a,b}f_a^\prime f_b^\prime$.

\medbreak\noindent We compute:
$$
\begin{array}{rl}
\psi_1^0=&1,\\
\psi_1^1=&\{\nabla_{\partial r}(r^\alpha\psi_1)\}|_{\partial M}=\{(\partial_r
-\textstyle{\textstyle\frac12}\sum_af_a^\prime)(1)\}|_{\partial M}=
-{\textstyle\frac12}\sum_af_a^\prime,\vphantom{\vrule height 11pt}\\
\textstyle\psi_1^2=&{\textstyle\frac12}\{(\nabla_{\partial_r})^2(r^\alpha\psi_1)\}|_{\partial M}
={\textstyle\frac12}\{(\partial_
-\textstyle{\textstyle\frac12}\sum_af_a^\prime)^2(1)\}|_{\partial M}\vphantom{\vrule height 11pt}\\
\textstyle=&{\textstyle\frac18}\sum_{a,b}f_a^\prime f_b^\prime
-{\textstyle\frac14}\sum_af_a^{\prime\prime},\vphantom{\vrule
height 11pt}\\
\psi_2^0=&1,\vphantom{\vrule height 11pt}\\
\textstyle\psi_2^1=&\{\tilde\nabla_{\partial r}(\psi_2)\}|_{\partial M}
=\{(\partial_r+{\textstyle\frac12}\sum_af_a^\prime)(e^{-\sum_af_a})\}|_{\partial
M}=-{\textstyle\frac12}\sum_af_a^\prime,\vphantom{\vrule height 11pt}\\
\textstyle\psi_2^2=&{\textstyle\frac12}\{(\tilde\nabla_{\partial r})^2\psi_2\}|_{\partial M}
={\textstyle\frac12}\{(\partial_r
+{\textstyle\frac12}\sum_af_a^\prime)^2(e^{-\sum_af_a})\}|_{\partial M}\vphantom{\vrule height 11pt}\\
=&\textstyle{\textstyle\frac18}\sum_{a,b}f_a^\prime f_b^\prime
-{\textstyle\frac14}\sum_af_a^{\prime\prime}\,.\vphantom{\vrule
height 11pt}
\end{array}$$

Considering the term $\sum_af_a^\prime$ in $\beta_{1,\alpha}^{\partial M}$ yields Assertion (b),
considering the term $\sum_a\delta_a^2$ in
$\beta_{2,\alpha}^{\partial M}$ yields Assertion (c), considering the term $\sum_af_a^{\prime\prime}$  in
$\beta_{2,\alpha}^{\partial M}$ yields Assertion (d), considering the term $\sum_{a,b}f_a^\prime f_b^\prime$
in $\beta_{2,\alpha}^{\partial M}$ yields Assertion (e), and considering the term $\sum_a(f_a^\prime)^2$ in
$\beta_{2,\alpha}^{\partial M}$ yields Assertion (f).

\subsection{The Proof of Theorem~\ref{thm-3}} We must now simply trace through the logic train.  We have
computed that:
\medbreak\quad
$\varepsilon_{\alpha_1,\alpha_2}=c_{\alpha_1,\alpha_2}$,
\medbreak\quad
$\varepsilon_{\alpha_1,\alpha_2}^1=c_{\alpha_1-1,\alpha_2}$ and $\varepsilon_{\alpha_1,\alpha_2}^3
=c_{\alpha_1,\alpha_2-1}$,
\medbreak\quad
$\varepsilon_{\alpha_2,\alpha_2}^2=-\frac12(\varepsilon_{\alpha_1,\alpha_2}^1
+\varepsilon_{\alpha_1,\alpha_2}^3)
 =-\frac12(c_{\alpha_1-1,\alpha_2}+c_{\alpha_1,\alpha_2-1})$,
\medbreak\quad
$\varepsilon_{\alpha_1,\alpha_2}^4=c_{\alpha_1-2,\alpha_2}$ and
$\varepsilon_{\alpha_1,\alpha_2}^7=c_{\alpha_1,\alpha_2-2}$,
\medbreak\quad
$\varepsilon_{\alpha_1,\alpha_2}^6=c_{\alpha_1,\alpha_2}$ and
$\varepsilon_{\alpha_1,\alpha_2}^{14}=c_{\alpha_1-1,\alpha_2-1}$,
\medbreak\quad
$\varepsilon_{\alpha_1,\alpha_2}^{12}=-\varepsilon_{\alpha_1,\alpha_2}^6=-c_{\alpha_1,\alpha_2}$,
\medbreak\quad
$\varepsilon_{\alpha_1,\alpha_2}^5=\varepsilon_{\alpha_1-1,\alpha_2}^2=
-\frac12(c_{\alpha_1-2,\alpha_2}+c_{\alpha_1-1,\alpha_2-1})$,
\medbreak\quad
$\varepsilon_{\alpha_1,\alpha_2}^8=\varepsilon_{\alpha_1,\alpha_2-1}^2=
-\frac12(c_{\alpha_1-1,\alpha_2-1}+c_{\alpha_1,\alpha_2-2})$,
\medbreak\quad
$\varepsilon_{\alpha_1,\alpha_2}^{11}=\varepsilon_{\alpha_1,\alpha_2}^9
=-\frac14\varepsilon_{\alpha_1,\alpha_2}^4-\frac14\varepsilon_{\alpha_1,\alpha_2}^7
+\frac12\varepsilon_{\alpha_1,\alpha_2}^6$
\smallbreak\qquad\qquad
$=-\frac14c_{\alpha_1-2,\alpha_2}-\frac14c_{\alpha_1,\alpha_2-2}+\frac12c_{\alpha_1,\alpha_2}$,
\medbreak\quad
$\varepsilon_{\alpha_1,\alpha_2}^{10}=-\{\textstyle\frac18\varepsilon_{\alpha_1,\alpha_2}^4
+{\textstyle\frac12}\varepsilon_{\alpha_1,\alpha_2}^5+{\textstyle\frac14}\varepsilon_{\alpha_1,\alpha_2}^6
+{\textstyle\frac18}\varepsilon_{\alpha_1,\alpha_2}^7
+{\textstyle\frac12}\varepsilon_{\alpha_1,\alpha_2}^{8}
+\frac14\varepsilon_{\alpha_1,\alpha_2}^{14}\}$
\medbreak\qquad\qquad
$=-\{\frac18c_{\alpha_1-2,\alpha_2}-\frac14(c_{\alpha_1-2,\alpha_2}+c_{\alpha_1-1,\alpha_2-1})
+\frac14c_{\alpha_1,\alpha_2}$
\medbreak\qquad\qquad\qquad
$+\frac18c_{\alpha_1,\alpha_2-2}-\frac14(c_{\alpha_1-1,\alpha_2-1}+c_{\alpha_1,\alpha_2-2})
+\frac14c_{\alpha_1-1,\alpha_2-1}\}$
\medbreak\qquad\qquad
$=\{(-\frac18+\frac14)c_{\alpha_1-2,\alpha_2}+(\frac14-\frac14+\frac14)c_{\alpha_1-1,\alpha_2-1}
+(-\frac18+\frac14)c_{\alpha_1,\alpha_2-2}$
\medbreak\qquad\qquad\qquad
$-\frac14c_{\alpha_1,\alpha_2}\}$.
\hfill\qedbox\medbreak

\section*{Acknowledgements} Research partially supported by project MTM2009-07756
(Spain) and by project 174012 (Serbia).

\end{document}